\documentclass{amsart}
\usepackage{amsmath}
\usepackage{graphics, epsfig}

\theoremstyle{plain}
\newtheorem{thm}{Theorem}[section]
\newtheorem{lem}[thm]{Lemma}

\theoremstyle{definition}
\newtheorem{defn}[thm]{Definition}

\def \CPb {\overline{\mathbb{CP}}^{2}}
\def \CP {{\mathbb{CP}}^{2}} 
\def \R {\mathbf{R}}
\def \Z {\mathbf{Z}}
\def \Sig{\Sigma}

\def \b {\beta}

\def \s {\sigma}

\def \bd {\partial}

\def \- {\setminus}

\def \ssw {\text{SW}}

\def\fs{\mathfrak{s}}

\title{Small exotic 4-manifolds} 
              
\begin{document}

\author{Anar Akhmedov}
\address{School of Mathematics\\
Georgia Institute of Technology\\
Atlanta, GA 30322}

\email{ahmadov@math.gatech.edu}

\dedicatory{Dedicated to Professor Ronald J. Stern on the occasion of his sixtieth birthday}

\begin{abstract} In this article, we construct the first examples of a simply connected minimal symplectic 4-manifold homeomorphic but not diffeomorphic to $3\CP\,\#7\CPb$. We also construct the first exotic \emph{symplectic} structure on $\CP\#5\,\CPb$.    
\end{abstract}

\maketitle

\setcounter{section}{-1}

\section{Introduction}

\medskip

Over the past several years, there has been a considerable progress in the discovery of exotic smooth structures on simply-connected $4$-manifolds with small Euler characteristic. In 2004, Jongil Park \cite{P2} has produced the first example of exotic smooth structure on $\CP\#7\,\CPb$ i.e. 4-manifold homeomorphic to $\CP\#\, 7\CPb$ but not diffeomophic to it. Soon after, Andr\'as Stipsicz and Zolt\'an Szab\'o used a similar technique to construct an exotic smooth structure on  $\CP\#\, 6\CPb$ \cite{SS1}. Then Fintushel and Stern \cite{FS3} introduced a new technique, surgery in double nodes, which demonstrated that in fact $\CP\#k\,\CPb$, $k= 6,7,8$, have infinitely many distinct smooth structures. Park, Stipsicz, and Szab\'o  \cite{PSS}, using double node surgery technique \cite{FS3}, constructed infinitely many smooth structures when $k = 5$. Stipsicz and Szab\'o used similar  ideas to construct the exotic smooth structures on $3\,\CP\#k\,\CPb$ for $k = 9$ \cite{SS2} and Park for $k = 8$ \cite{P3}. In this article, we construct an exotic smooth structure on $3\,\CP\#7\,\CPb$. We also construct an exotic \emph{symplectic} structure on $\CP\#5\,\CPb$, the first known such symplectic example.  

Our approach is different from the recent constructions in the sense that we do not use any rational-blowdown surgery (\cite{FS4}, \cite{P1}). The techniques used in our construction are the symplectic fiber sum operation (\cite{G}, \cite{MW}) and the symplectic cohomology $S^{2}\times S^{2}$ \cite{A1}, the recent construction of the author. Our results are the following

\medskip

\begin{thm}

There exist a smooth closed simply-connected minimal symplectic 4-manifold $X$ that is homeomorphic but not diffeomorphic to $3\,\CP\#7\,\CPb$.

\end{thm}

\medskip

\begin{thm}

There exist a smooth closed simply-connected minimal symplectic 4-manifold $Y$ which is homeomorphic, but not diffeomorphic, to a rational surface $\CP\#5\,\CPb$.

\end{thm}

The article is organized as follows. The first two sections give a quick introduction to Seiberg-Witten invariants and a fiber sum operation. In the third section, we review the building blocks for our construction. In sections 4 and 5 we construct symplectic 4-manifolds $X$ and $Y$ homeomorphic but not diffemorphic $3\,\CP\#7\,\CPb$ and $\CP\#5\,\CPb$, respectively.    

\medskip

\noindent {\bf Acknowledgments:} I would like to thank John Etnyre, Ron Stern and Andr\'as Stipsicz for their interest in this work, and for their encouragement. Also, I am very grateful to B. Doug Park for the comments on the first draft of this article, kindly pointing out some errors in the fundamental group computations and for his corrections. This work is partially supported by NSF grant FRG-0244663.

\section{Seiberg-Witten Invariants}

 In this section we review the basics of Seiberg-Witten invariants introduced by Seiberg and Witten. Let us recall that the Seiberg-Witten invariant  of a smooth closed oriented $4$-manifold $X$ with $b_2 ^+(X)>1$ is an integer valued function which is defined on the set of $spin ^{\, c}$ structures over $X$ \cite{W}. For simplicity we assume that $H_1(X,\Z)$ has no 2-torsion. Then there is a one-to-one correspondence between the set of $spin ^{\, c}$ structures over $X$ and the set of characteristic elements of $H^2(X,\Z)$ as follows: to each $spin ^c$ structure $\fs $ over $X$ corresponds a bundle of positive spinors $W^+_{\fs}$ over $X$. Let $c(\fs)\in H_2(X, \Z)$ denote the Poincar\'e dual of $c_1(W^+_{\fs})$. Each of $c(\fs)$ is a characteristic element of $H_2(X,\Z)$ (i.e. its Poincar\'e dual $\hat{c}(\fs)=c_1(W^+_{\fs})$ reduces mod~2 to $w_2(X)$).

In this set up we can view the Seiberg-Witten invariant as integer valued function
\[ \ssw_X: \lbrace k\in H^2(X,\Z)|k\equiv w_2(TX)\pmod2 \rbrace
\rightarrow \Z. \] The Seiberg-Witten invariant $\ssw_X$ is a
diffeomorphism invariant and its sign depends on an orientation of
\[ H^0(X,\R)\otimes\det H_+^2(X,\R)\otimes \det H^1(X,\R).\] If
$\ssw_X(\b)\neq 0$, then we call
$\b$ a {\it{basic class}} of $X$. It is a fundamental fact that the set
of basic classes is finite. It can be shown that, if $\b$ is a basic class, then
so is $-\b$ with
\[\ssw_X(-\b)=(-1)^{(e+\s)(X)/4}\,\ssw_X(\b)\] where
$e(X)$ is the Euler characteristic and $\sigma(X)$ is the signature
of $X$.

\medskip

\begin{thm}\cite{T} Suppose that $X$ is a closed symplectic 4-manifold with $b_{2}^{+}(X) > 1$. If $K_{X}$ is a canonical class of $X$, then $\ssw_{X}(\pm K_{X}) = \pm 1$. \end{thm}

\medskip

\begin{thm}(Liu, Ohta-Ono). Let $X$ be a closed minimal symplectic 4-manifold with $b_{2}^{+}(X) = 1$ and a canonical class $K_{X}$. Then the followings are equivalent

(i) $X$ admits a metric of positive scalar curvature.

(ii) $X$ admits a symplectic structure $\omega$ with $K_{X}\cdot [\omega] < 0$

(iii) $X$ is either rational or ruled.

\end{thm}

\section{Fiber Sum}

\begin{defn} Let $X$\/ and $Y$\/ be closed, oriented, smooth $4$-manifolds each containing a smoothly embedded surface $\Sigma$ of genus $g \geq 1$. Assume $\Sig$ represents a homology of infinite order and has self-intersection zero in $X$\/ and $Y$\/, so that there exist a tubular neighborhood, say $\nu\Sigma\cong \Sigma\times D^{2}$, in both $X$\/ and $Y$. Using an orientation-reversing and fiber-preserving diffeomorphism $ \psi : S^{1}\times\Sigma  \longrightarrow   S^{1}\times\Sigma$, we can glue $X \setminus\nu\Sigma$ and $Y\setminus \nu\Sigma$\/ along the boundary $\partial(\nu\Sigma)\cong \Sigma\times S^{1}$. This new oriented smooth $4$-manifold $X\#_{\psi}Y$ is called a \emph{generalized fiber sum}\/ of $X$\/ and $Y$\/ along $\Sigma$.  

\end{defn}

\medskip

\begin{lem} 
Let $X$ and $Y$ be closed, oriented, smooth\/ $4$-manifolds containing an embedded surface\/ $\Sigma$ of self-intersection\/ $0$. Then 
\begin{eqnarray*}
c_{1}^{2}(X\#_{\psi}Y) &=&  c_{1}^{2}(X) + c_{1}^{2}(Y) + 8(g-1),\\
\chi_{h}(X\#_{\psi}Y) &=&  \chi_{h}(X) + \chi_{h}(Y) + (g-1),
\end{eqnarray*}
where $g$ is the genus of the surface $\Sigma$. 
\end{lem}

\begin{proof}
The above formulas simply follow from the well-known formulas 
\begin{equation*}
e(X\#_{\psi}Y)= e(X) + e(Y) - 2e(\Sigma),\quad 
\sigma(X\#_{\psi}Y) = \sigma(X) + \sigma(Y).  \qedhere   
\end{equation*} once we apply the formulas $\chi_{h} = (\sigma  + e) / 4$ and  ${c_{1}^{2}}= 3\sigma + 2e$
\end{proof}
 
\noindent If $X$, $Y$ are symplectic manifolds and $\Sig$ is a symplectic submanifold then according to theorem of Gompf \cite{G} $X\#_{\psi}Y$ admits a symplectic structure.  

\noindent We will use the following recent theorem  of M. Usher \cite{U} to show that the symplectic manifolds constructed in Sections 4 and 5 are minimal. 

\medskip

\begin{thm}\cite{U} {\bf (Minimality of Sympletic Sums)} Let $Z = X_{1}\#_{F_{1} = F_{2}}X_{2}$ be sympletic fiber sum of manifolds $X_{1}$ and $X_{2}$. Then:

(i) If either $X_{1} \backslash F_{1}$ or $X_{2} \backslash F_{2}$ contains an embedded sympletic sphere of square $-1$, then $Z$ is not minimal.

(ii) If one of the summands $X_{i}$ (say $X_{1}$) admits the structure of an $S^{2}$-bundle over a surface of genus $g$ such that $F_{i}$ is a section of this fiber bundle, then $Z$ is minimal if and only if $X_{2}$ is minimal. 

(iii) In all other cases, $Z$ is minimal.

\end{thm}

\section {Main Building blocks}

Our building blocks will be the symplectic cohomology $S^{2}\times S^{2}$ \cite {A1}, the recent construction of the author and the manifold $T^{2}\times S^{2}\,\#4\CPb$. 

\medskip

\subsection{Matsumoto fibration}

First, recall that the manifold $Z = T^{2}\times S^{2}\,\#4\CPb$ can be described as the double branched cover of $S^{2}\times T^{2}$ where branch set $B_{2, 2}$ is the union of two disjoint copies of $S^{2}\times \{{\rm{pt}}\}$ and two disjoint copies of $\{{\rm{pt}}\}\times T^{2}$. The branch cover has $4$ singular points corresponding to number of intersection of horizontal lines and vertical tori in the branch set $B_{2, 2}$. After desingularizing the above singular manifold, one obtains $T^{2}\times S^{2}\,\#4\CPb$. The vertical fibrations of $S^{2}\times T^{2}$ pull back to give fibrations of $T^{2}\times S^{2}\,\#4\CPb$ over $S^{2}$. A generic fiber of the vertical fibration is the double cover of $T^2$, branched over $2$ points. Thus a generic fiber will be a genus two surface. According to Matsumoto \cite{M} , this fibration can be perturbed to be a Lefschetz fibration over $S^{2}$ with global monodromy $(\beta_{1} \beta_{2} \beta_{3} \beta_{4})^{2} = 1$, where the curves $\beta_{1}$, $\beta_{2}$, $\beta_{3}$ and $\beta_{4}$ are shown in Figure 1.

\begin{figure}[htpb]

\centering 
\psfig{file=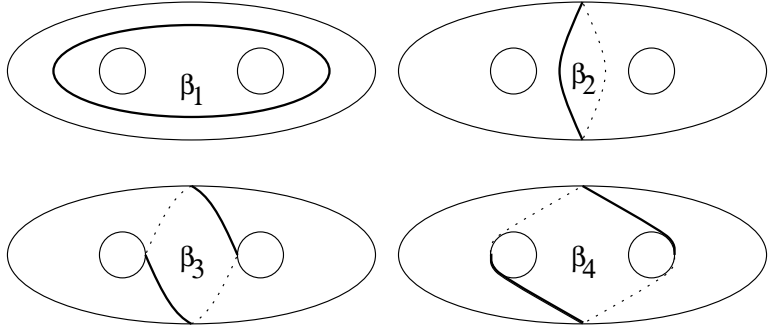, width = 4.0in, clip=}
\caption{Dehn Twists for Matsumoto fibration}
\label{matsumotocurves}
\end{figure}

Let us denote a regular fiber as $\Sig'_{2}$ and the standard generators of fundamental group of $\Sig'_{2}$ as ${a_{1}}$, ${b_{1}}$, ${a_{2}}$ and ${b_{2}}$. Using the homotopy exact sequence for a Lefschetz fibration 

\begin {center} 

$\pi_{1}(\Sig'_{2}) \longrightarrow \pi_{1}(Z) \longrightarrow \pi_{1}(S^{2})$

\end{center}

\noindent According to \cite{OS}, we have the following identification of the fundamental group of $Z$: 

\medskip

$\pi_{1}(Z) = \pi_{1}(\Sig'_{2})/<\beta_{1},\beta_{2},\beta_{3},\beta_{4}> $ 

\medskip

\noindent $\beta_{1} = b_{1}b_{2}$ \

\noindent $\beta_{2} = a_{1}b_{1}{a_{1}}^{-1}{b_{1}}^{-1} = a_{2}b_{2}{a_{2}}^{-1}{b_{2}}^{-1}$ \

\noindent $\beta_{3} = b_{2}a_{2}{b_{2}}^{-1}a_{1}$ \

\noindent $\beta_{4} = b_{2}a_{2}a_{1}b_{1}$ \

\medskip

\noindent Hence $\pi_{1}(Z) = \ < a_{1}, b_{1}, a_{2}, b_{2} \ | \  b_{1}b_{2} = [a_{1}, b_{1}] = [a_{2}, b_{2}] = b_{2}a_{2}{b_{2}}^{-1}a_{1} = 1 >$

\medskip

Note that the fundamental group of $T^{2}\times S^{2}\,\#4\CPb$ is ${\bf Z\oplus Z}$, generated by two of these generators (say $a_{1}$ and $b_{1}$). The other two generators $a_{2}$ and $b_{2}$ are the inverses of $a_{1}$ and $b_{1}$ in the fundamental group. Also, the fundamental group of the complement of $\nu \Sig'_{2}$ is ${\bf Z\oplus Z}$. It is generated by $a_{1}$ and $b_{1}$. The normal circle $\lambda' = {pt} \times \bd D^{2}$ to $\Sig'_{2}$ can be deformed using an exceptional sphere section, thus trivial in $\pi_{1}(T^{2}\times S^{2}\,\#4\CPb \setminus \nu \Sig'_{2})$.

\medskip  

\begin{lem} ${c_1}^{2}(Z) = -4$, $\chi_h(Z) = 0$ and $K_{Z} = 2U + E_{1} + E_{2} + E_{3} + E_{4}$, where $K_{Z}$ is canonical class of the symplectic structure on $Z$, $U$ is the class of $T^{2}$ and $E_{i}$ are the exceptional classes. 

\end{lem}

\begin{proof}

We have ${c_1}^{2}(Z) = {c_1}^{2}(T^{2} \times S^{2}) - 4 = -4$, $\sigma (Z) = \sigma(T^{2}\times S^{2}) - 4 = -4$ and $\chi_{h}(Z) = \chi_{h}(T^{2}\times S^{2}) = 0$. The latter follows from the generalized adjunction inequality and the blowup formula for Seiberg-Witten invariants. \end{proof}

\subsection {Symplectic 4-manifolds cohomology equivalent to $S^{2}\times S^{2}$}

Our second building block will be the symplectic cohomology $S^{2}\times S^{2}$ \cite{A1}. For the sake of completeness, the details of this construction are included below. We refer the reader to \cite{A1} for more details.

Let $K$ be a genus one fibered knot (i.e. trefoil or figure eight knot) in $S^{3}$ and $m$ a meridional circle to $K$. We perform 0-framed surgery on $K$ and denote the resulting 3-manifold by $M_{K}$. The manifold $M_{K}$ has same integral homology as $S^{2}\times S^{1}$, where the class of $m$ generates $H_{1}$. Since the knot is genus $1$ fibered knot, it follows that the manifold $M_{K}\times S^{1}$ is a torus bundle over the torus which is homology equivalent to $T^{2}\times S^{2}$. Since $K$ is fibered knot, $M_{K}\times S^{1}$ admits a symplectic structure. Note that there is section $\ s : S^{1}\times S^{1} \longrightarrow M_{K}\times S^{1}$ of this fibration. Both the torus fiber and the torus section are symplectically embedded and have a self-intersection zero. The first homology of $M_{K}\times S^{1}$ is generated by the standard first  homology generators $m$ and $x$ of the torus section. The classes of circles $\gamma_{1}$ and $\gamma_{2}$ of fiber $F$ of the given fibration are trivial in homology.

We form a twisted fiber sum of the two copies of the above manifold $M_{K}\times S^{1}$  where we identify the fiber of one fibration to the base of the other. Let $Y_{K}$ denote the mentioned twisted fiber sum $Y_{K} = M_{K}\times S^{1} \#_{F = s(T_{m})}\,\ M_{K}\times S^{1}$. It follows from Gompf's theorem \cite {G} that $Y_{K}$ is symplectic. 

Let $T_{1}$ be the section of the first copy of $M_{K}\times S^{1}$ in $Y_{K}$ and $T_{2}$ be the fiber of the second copy. Then $T_{1} \# T_{2}$ symplectically embeds in $Y_{K}$ \cite {FS2}. Now suppose that $Y_{K}$ is the symplectic 4-manifold given, and $\Sig _{2} = T_{1} \# T_{2}$ is the genus two symplectic submanifold of self-intersection zero sitting inside of $Y_{K}$. Let $m$, $x$, $\gamma_{1}$ and $\gamma_{2}$ be the generators of $H_{1}(\Sig_{2}, \Z)$. We choose the diffeomorphism  $\phi : T_{1} \# T_{2} \longrightarrow T_{1} \# T_{2}$ of the $\Sig_{2}$ that changes the generators of the first homology according to the following rule: $\phi (m') = \gamma_{1}$, $\phi ({\gamma_1}') = m$, $\phi (x') = \gamma_{2}$ and $\phi ({\gamma_2}') = x$. Next we take the fiber sum of two copies of $Y_{K}$ along the genus two surface $\Sig_{2} = S$ using the diffeomorphism $\phi$. Let denote the new symplectic manifold as $X_{K}$ i.e. $X_{K} = Y_{K}\#_{\Sig_{2} = \phi (\Sig_{2})}\,\ Y_{K}$. In \cite{A1} we show that the manifold $X_{K}$ has a trivial first Betti number and has integral cohomology of $S^{2}\times S^{2}$. $H_{2}(X_{K},{\Z}) = \Z \oplus \Z$, where base for the second homology are classes of the $\Sig_{2} = S$ and the new genus two surface $T$ resulting from the fiber sum. Also, $S^{2} = T^{2} = 0$ and $S \cdot T = 1$. Furthermore ${c_1}^{2}(X_{K}) = 8$, $\sigma (X_{K}) = 0$ and $\chi_{h}(X_{K}) = 1$. It follows from Theorem 2.3 that $Y_{K}$ and $X_{K}$ are both minimal. 

\subsubsection{Step 1: Fundamental Group of $M_{K}\times S^{1}$} We will assume that $K$ is trefoil. Let $a$, $b$ and $c$ denote the Wirtinger generators of the trefoil knot. The knot group of the trefoil has the following presentation $ < a, b, c  \ | \  abc^{-1}{b}^{-1}, \ ca{b}^{-1}{a}^{-1}>$ \hspace{1mm}  =  \hspace{1mm} $< a, b  \ | \  aba{b}^{-1}{a}^{-1}{b}^{-1}>$ \hspace{1mm} = \hspace{1mm}  $ <u, v  \ | \  u^{2} =  v^{3} >$ where $u = bab$ and $v = ab$. The homotopy classes of the meridian and the longitude of the trefoil are given as follows: $ m = uv^{-1} = b$ and $l = u^{2}(uv^{-1})^{-6} = ab^2ab^{-4}$ \cite{BZ}. Also, the homotopy classes of $\gamma_{1} = {a}^{-1}{b}$ and $\gamma_{2}= {b}^{-1}{a}ba^{-1}$. Notice that the fundamental group of $M_{K}$, 0-surgery on the trefoil, is obtained from the knot group of the trefoil by adjoining the relation $l = u^{2}(uv^{-1})^{-6} = ab^2ab^{-4} = 1$, i.e. $\pi_{1}(M_{K}) \ = \ <u, v  \ | \  u^{2} = v^{3}, \ u^{2}(uv^{-1})^{-6} = 1 > \ = \ < a, b \ | \ aba = bab, \ ab^{2}a = b^4 >$. $\pi_{1}(M_{K}\times S^{1}) \hspace{1mm} = \hspace{1mm}  < a, b, x  \ | \  aba = bab, \ ab^{2}a = b^4, \ [x, a] = [x, b] = 1 >$.         

\subsubsection{Step 2: Fundamental Group of $Y_{K}$} Next we take two copies of the manifold $M_{K}\times S^{1}$. In the first copy, take a tubular neighborhood of the torus section $T_{m}$, remove it from $M_{K}\times S^{1}$ and denote the resulting manifold as $C_{B}$. In the second copy, we remove a tubular neighborhood of the fiber $F$ and denote it by $C_{F}$. Let $Y_{K}$ denote the fiber sum $Y_{K} = M_{K}\times S^{1} \#_{F\times D^{2} = s(T_{m})\times D^{2}}\,\ M_{K}\times S^{1}$, i.e. we glue $C_{B}$ and $C_{F}$ along their common boundary $T^{3}$. Notice that $C_{B} = {M_{K} \times S^{1}} \backslash {T_{m} \times D^{2}} =  {(M_{K} \backslash S^{1} \times D^{2})} \times S^{1}$. We have $\pi_{1}(C_{B}) = \pi_{1}(K) \ \oplus \ <x>$ where $x$ is the generator corresponding to the $S^{1}$ copy. The new circle $\lambda$, resulting from the removal of the torus $S^{1} \times D^{2} = m\times D^{2}$ from the 3-manifold $M_{K}$, is the longitude of the knot $K$, thus trivial in the homology. $H_{1}(C_{B}) = H_{1}(M_{K}) = \ <m> \oplus <x>$. 

To compute the fundamental group of the $C_{F}$, we will use the following observation: $F\times D^2 = s^{-1}(D^{2})$ i.e. it is the preimage of the small disk on $T_{m'} = m' \times y$ and $T_{m'} \backslash D^{2}$ is homotopy equivalent to the wedge of two circles. The generators of the fundamental group $y$ and $m'$ of the new base do not commute anymore, but $y$ still commutes with generators $\gamma_{1}'$ and $\gamma_{2}'$ of the $\pi_{1}(F)$ which maps into $\pi_{1}(M_{K}\times S^{1})$. The fundamental group and the first homology of the $C_{F}$ will be isomorphic to the followings: $ \pi_{1}(C_{F}) = \ < d, y, \gamma_{1}', \gamma_{2}' \ | \ [y, \gamma_{1}'] = [y, \gamma_{2}'] = [ \gamma_{1}', \gamma_{2}'] = 1, \ d {\gamma_{1}'}d^{-1} = \gamma_{1}' \gamma_{2}', \ d{\gamma_{2}'}d^{-1} = ({\gamma_{1}'})^{-1}>$ and $H_{1}(C_{F}) = \ <d> \oplus <y>$.            

We use the Van Kampen's Theorem to compute the fundamental group of $Y_{K}$

\noindent $ \pi_{1}(Y_{K}) \ = \pi_{1}(C_{F})_ {*_{\pi_{1}(T^3)}} \pi_{1}(C_{B}) \ =  \ < d, y, \gamma_{1}', \gamma_{2}' \ | \ [y, \gamma_{1}'] = [y, \gamma_{2}'] = [ \gamma_{1}', \gamma_{2}'] = 1, \ d {\gamma_{1}'}d^{-1} = \gamma_{1}' \gamma_{2}', \ d{\gamma_{2}'}d^{-1} = ({\gamma_{1}}')^{-1} > _{< \gamma_{1}' = x,\  \gamma_{2}' = b, \ {\lambda}' = {\lambda} >} < a, b, x  \ | \ aba =bab,  \ [x, a] =  [x, b] = 1 > \ = \ < a, b, x, d, y \ | \  aba = bab, \ [x, a] = [x, b] = 1, \ [y, x] = [y, b] = 1, \ dxd^{-1} = xb, \ dbd^{-1} = x^{-1}, \ [d, y] = ab^{2}ab^{-4}> $.

\subsubsection {Step 3: Fundamental Group of $X_{K}$} 

Finally, we carry out the computations of the fundamental group and the the first homology of $X_{K}$. Let $Y_{1} = Y_{2} = Y_{K} \setminus (T_{1}\# T_{2})\times D^{2}$. Again, by Van Kampen's theorem we have \

$ \pi_{1}(X_{K}) = \pi_{1}(Y_{1}) *_{ (T_{1}\#T_{2})\times D^2 = \phi (T_{1}\#T_{2}\times D^{2})} \pi_{1}(Y_{2}) \ = $  \\			
$ < a, b, x, d, y; e, f, z, s, t; \ l_{1}, \cdots , l_{m}; \ {l_1}',\cdots , {l_m}' \ | \  aba =bab, \ [y, x] = [y, b] = 1, \ dxd^{-1} = xb, \ dbd^{-1} = x^{-1}, \ ab^{2}ab^{-4} = [d, y],\ r_{1} = \cdots = r_{n} = 1, \ {r_1}' = \cdots = {r_n}'= 1, \ efe = fef, \ [t, z] = [t, f] = 1, \ szs^{-1} = zf, \ sfs^{-1} = z^{-1}, \ ef^{2}ef^{-4} = [s, t], \ d = e^{-1}f, \ y = f^{-1}efe^{-1}, \ a^{-1}b = s, \ b^{-1}aba^{-1} = t, \ [x, b] = [z, f] >$ 

\noindent where the elements $l_{i}, {l_i}'$ (for $i = 1, \cdots , m $) and $r_{j}, {r_{j}}'$ (for $j = 1, \cdots , n$) all are in the normal subgroup generated by $[x, b] = [z, f]$
 
\medskip

 By abelianizing $\pi_{1}(X_{K})$, we have $H_{1}(X_{K}, \Z ) = 0$

Notice that it follows from our gluing that the images of standard generators of the fundamental group of $\Sig_{2}$ are  $a^{-1}b$, $b^{-1}aba^{-1}$, $d$ and $y$ in $\pi_{1}(X_{K})$. 

\medskip

\begin{lem} $K_{Y_{K}} = 2F$, where $F$ is the class of fiber in $Y_{K}$ and $K_{Y_{K}}$ is canonical class of the symplectic structure on $Y_{K}$. 

\end{lem}

\begin{proof} Notice that $K_{M_{K}\times S^{1}} = 0$ when $K$ is trefoil. By the canonical class formula for a fiber sum, we have the following

$K_{Y_K} = 2K_{M_{K} \times S^{1}} + 2F = 2F$ \end{proof}

\medskip

\section{Construction of exotic $3\,\CP\#7\,\CPb$}

In this section, we construct a simply-connected minimal symplectic $4$-manifold $X$ homeomorphic but not diffemorphic to $3\,\CP\#7\,\CPb$. Using Seiberg-Witten invariants, we will distinguish $X$ from $3\,\CP\#7\,\CPb$. 

Our manifold $X$ will be the symplectic fiber sum of $X_{K}$ and $Z = T^{2}\times S^{2}\,\#4\CPb$ along the genus two surfaces $\Sigma_{2} = S$ and $\Sigma'_{2}$. Recall from \cite{A1} that  $a^{-1}b$, $b^{-1}aba^{-1}$, $d$, $y$, and 
$\lambda=\{{\rm pt}\}\times S^1=[x,b][z,f]^{-1}$ generate the 
inclusion-induced image of $\pi_1(\Sigma_{2}\times S^{1})$ 
inside $\pi_1(X_K\setminus\nu\Sigma_2)$.  
As before, let $a_{1}$, $b_{1}$, $a_{2}$, $b_{2}$ and $\lambda'$ generate $\pi_1(\Sigma'_{2}\times S^{1})$. We choose the gluing diffeomorphism  $\psi : \Sigma_{2}\times S^1 \rightarrow \Sigma'_{2}\times S^1$ that maps the fundamental group generators as follows: 
\begin{equation*}
\psi_{\ast} (a^{-1}b) = a_2,\;
\psi_{\ast} (b^{-1}aba^{-1}) = b_2,\;
\psi_{\ast} (d) = a_1 ,\; \psi_{\ast} (y) = b_1 ,\;  \psi_{\ast} (\lambda) = \lambda' .
\end{equation*} 
It follows from Gompf's theorem \cite{G} that $X = X_{K}\#_{\psi}(T^{2}\times S^{2}\,\#4\CPb)$ is symplectic. 
\medskip

\begin{lem} $X$ is simply connected. 

\end{lem}

\begin{proof}

By Van Kampen's theorem, we have \

\begin{equation*}
\pi_{1}(X) \:=\: \frac{\pi_{1}(X_{K}\setminus \nu \Sigma_{2})\ast \pi_{1}(Z \setminus \nu \Sigma'_{2})}{\langle a^{-1}b = a_{2},\, b^{-1}aba^{-1} = b_{2},\, d = a_{1},\, y = b_{1},\, 
\lambda = 1 \rangle}.
\end{equation*}

Note that a nontrivial element $\lambda$ of $\pi_{1}(X_{K} \setminus \nu \Sigma_{2})$ becomes trivial in $X$. Also, using the relations  $b_{1}b_{2} = [a_{1}, b_{1}] = [a_{2}, b_{2}] = b_{2}a_{2}{b_{2}}^{-1}a_{1} = a_{1}a_{2} = 1$, we get the following relations in the fundamental group of $X$: $a^{-1}bd = [a^{-1}b, b^{-1}aba^{-1}] = [d, y] = [d, b^{-1}aba^{-1}] = yb^{-1}aba^{-1} = 1$. Notice that the fundamental group of $Z$ is an abelian group of rank two. In addition, we have the following relations in $\pi_{1}(X)$ coming from the fundamental group of $X_{K}$: $aba = bab$, $efe = fef$, $[y, b] = [t, f] = 1$, $dbd^{-1} = x^{-1}$, $dxd^{-1} = xb$, $sfs^{-1} = z^{-1}$, $szs^{-1} = zf$, $a^{-1}b = s$,  $b^{-1}aba^{-1} = t$, $y = f^{-1}efe^{-1}$ and $e^{-1}f = d$. These set of relations give rise to the following identities

\smallskip

\begin{eqnarray}
a &=& bd,\, \label{id:abd}\\
yb &=& by,\, \label{id:yby}\\
aba &=& bab,\, \label{id:aba}\\
yab &=& ba,\, \label{id:yab}
\end{eqnarray}

\medskip

\noindent Next multiply the relation $(4)$ by $a$ from the right and use $aba = bab$. We have $yaba = ba^{2}$ $\Longrightarrow$ $ybab = ba^{2}$. By cancelling the element $b$, $yab = a^{2}$. Finally applying the relation $(4)$ again, we have $ba = a^{2}$. Later implies that $b = a$. Since $a = bd$, $dbd^{-1} = x^{-1}$, $dxd^{-1} = xb$, $aba = bab$ and $yb^{-1}aba^{-1} = 1$,  we obtain $d = y = x = b = a = 1$. Furthemore, using the relations $a^{-1}b = s$,  $b^{-1}aba^{-1} = t$, $efe = fef$, $e^{-1}f = d$, $sfs^{-1} = z^{-1}$ and $szs^{-1} = zf$, we similarly have $s = t = z = f = e = 1$. Thus, we can conclude that the elements $a$, $b$, $x$ $d$, $y$, $e$, $f$, $z$, $s$ and $t$ are all trivial in the fundamental group of $X$. Since we identified $a^{-1}b$ and $b^{-1}aba^{-1}$ with generators $a_{2}$ and $b_{2}$ of the group $\pi_{1}(Z \setminus \nu \Sigma'_{2}) = {\bf Z\oplus Z}$, it follows that $a_{2}$ and $b_{2}$ are trivial in the fundamental group of $X$ as well. This proves that $X$ is simply connected.  

\end{proof}

\medskip

\begin{lem} ${c_1}^{2}(X) = 12$, $\sigma (X) = - 4$ and $\chi_h(X) = 2$

\end{lem}

\begin{proof} 

We have ${c_1}^{2}(X) = {c_1}^{2}(X_{K}) + {c_1}^{2}(T^{2}\times S^{2}\,\#4\CPb) + 8$, $\sigma (X) = \sigma(X_{K}) + \sigma(T^{2}\times S^{2}\,\#4\CPb)$ and $\chi_{h}(X) = \chi_{h}(X_{K}) + \chi_{h}(T^{2}\times S^{2}\,\#4\CPb) + 1$. Since ${c_1}^{2}(X_{K}) = 8$, $\sigma(X_{K}) = 0$ and $\chi_{h} (X_{K}) = 1$, the result follows from the lemmas 2.2 and 3.1. \end{proof}

By Freedman's theorem \cite{F} and the lemmas 4.1 and 4.2, we have $X$ is homeomorphic to $3\,\CP\#7\,\CPb$. It follows from Taubes theorem of Section 2 that $\ssw_{X}( K_{X} ) = \pm 1$. Next we apply the connected sum theorem for the Seiberg-Witten invariant and show that $SW$ function is trivial for $3\,\CP\#7\,\CPb$. Since the Seiberg-Witten invariants are diffeomorphism invariants, we conclude that $X$ is not diffeomorphic to $3\,\CP\#7\,\CPb$. Notice that case (i) of Theorem 2.3 does not apply and $X_{K}$ is minimal symplectic manifold. Thus, we can conclude that $X$ is minimal. Since symplectic minimality implies irreducibility for simply-connected $4$-manifolds with $b_{2}^{+} > 1$ \cite{K}, it follows that $X$ is also smoothly irreducible.    

\section{Construction of exotic symplectic $\CP\#5\,\CPb$}

In this section, we construct a simply-connected minimal symplectic $4$-manifold $Y$ homeomorphic but not diffemorphic to $\CP\#5\,\CPb$. Using Usher's Theorem \cite{U}, we will distinguish $Y$ from $\CP\#5\,\CPb$.   

The manifold $Y$ will be the symplectic fiber sum of $Y_{K}$ and $T^{2}\times S^{2}\,\#4\CPb$ along the genus two surfaces $\Sigma_{2} = S$ and $\Sigma'_{2}$. Let us choose the gluing diffeomorphism  $\varphi : \Sigma_{2}\times S^1\rightarrow \Sigma'_{2}\times S^1$ that maps the generators $a^{-1}b$, $b^{-1}aba^{-1}$, $d$, 
$y$\/ and $\mu$\/ of $\pi_{1}(Y_{K} \setminus \nu \Sigma_{2})$ to the generators 
$a_{1}$, $b_{1}$, $a_{2}$, $b_{2}$\/ and $\mu' = 1$ of 
$\pi_{1}(Z \setminus \nu \Sigma_{2}')$ according to the following rule: 
\begin{equation*}
\varphi_{\ast} (a^{-1}b) = a_2, \; \varphi_{\ast} (b^{-1}aba^{-1}) = b_2, \;
\varphi_{\ast} (d) = a_1, \; \varphi_{\ast} (y) = b_1, \;
\varphi_{\ast} (\mu) = \mu' = 1.
\end{equation*} 
Here, $\mu$ and $\mu'$ denote the meridians of $\Sigma$ and $\Sigma_2'$.

Again, by Gompf's theorem \cite{G}, $Y = Y_{K}\#_{\Sig}(T^{2}\times S^{2}\,\#4\CPb)$ is symplectic.  

\medskip

\begin{lem} $Y$ is simply connected. 

\end{lem}

\begin{proof}

By Van Kampen's theorem, we have \

\begin{equation*}
\pi_{1}(Y) \:=\: \frac{\pi_{1}(Y_{K}\setminus \nu \Sigma_{2})\ast \pi_{1}(Z \setminus \nu \Sigma'_{2})}{\langle a^{-1}b = a_{2},\, b^{-1}aba^{-1} = b_{2},\, d = a_{1},\, y = b_{1},\, 
\lambda = 1 \rangle}.
\end{equation*}

\noindent Using the exact same argument as in proof of Lemma 4.1, we have $a = b = x = d = y = 1$. Thus $\pi_{1}(Y) = 0$.   \end{proof}

\medskip

\begin{lem} ${c_1}^{2}(Y) = 4$, $\sigma(Y) = - 4$ and $\chi_h(Y) = 1$

\end{lem}

\begin{proof} 

We have ${c_1}^{2}(Y) = {c_1}^{2}(Y_{K}) + {c_1}^{2}(T^{2}\times S^{2}\,\#4\CPb) + 8$, $\sigma (Y) = \sigma(Y_{K}) + \sigma(T^{2}\times S^{2}\,\#4\CPb)$ and $\chi_{h}(Y) = \chi_{h}(Y_{K}) + \chi_{h}(T^{2}\times S^{2}\,\#4\CPb) + 1$. Since ${c_1}^{2}(Y_{K}) = 0$, $\sigma(Y_{K}) = 0$ and $\chi_{h} (Y_{K}) = 0$,  the result follows from lemmas 2.2 and 3.1. \end{proof}

By Freedman's classification theorem \cite{F} and lemmas 5.1 and 5.3 above, we have $Y$ is homeomorphic to $\CP\#5\,\CPb$. Notice that $Y$ is a fiber sum of non minimal manifold $Z = T^{2}\times S^{2}\,\#4\CPb$ with minimal manifold $Y_{K}$. All $4$ exceptional spheres $E_{1}$, $E_{2}$, $E_{3}$ and $E_{4}$ in $Z$ meet with the genus two fiber $2T + S - E_{1} - E_{2} - E_{3} - E_{4}$. It follows from Theorem 2.3. that $Y$ is a minimal symplectic manifold. Since symplectic minimality implies irreducibility for simply-connected $4$-manifolds for $b_{2}^{+} = 1$ \cite{HK}, it follows that $Y$ is also smoothly irreducible. We conclude that $Y$ is not diffeomorphic to $\CP\#5\,\CPb$.   

Theorem 1.2 provides an alternative way to show that $Y$ is exotic symplectic.

\end{document}